\documentclass{elsart}
\usepackage{amssymb}
\newtheorem{theo}{\arabic{section}.\arabic{abz}. Theorem}

\newtheorem{defi}{\arabic{section}.\arabic{abz}. Definition}

\newtheorem{nota}{\arabic{section}.\arabic{abz}. Notations}

\newtheorem{propo}{\arabic{section}.\arabic{abz}. Proposition}

\newtheorem{coro}{\arabic{section}.\arabic{abz}. Corollary}

\newtheorem{conv}{\arabic{section}.\arabic{abz}. Convention
and notations}

\newtheorem{lemm}{\arabic{section}.\arabic{abz}. Lemma}

\newtheorem{rema}{\arabic{section}.\arabic{abz}. Remark}

\newtheorem{assu}{\arabic{section}.\arabic{abz}.
Assumptions and notations}

\newtheorem{exa}{\arabic{section}.\arabic{abz}. Example}

\newcounter{abz}[section]
\newcounter{sabz}[abz]
\addtocounter{section}{-1}
\newcommand{\abz}{\refstepcounter{abz}}
\newcommand{\sabz}{\refstepcounter{sabz}}

\newcommand{\lbr}{\linebreak[0]}

\def\qed{Q.E.D.}
\def\rk{\mathop{\mathrm{rank}}}
\def\crk{\mathop{\mathrm{corank}}}

\def\im{\mathop{\mathrm{im}}}

\def\Ad{\mathop{\mathrm{Ad}}}
\def\ad{\mathop{\mathrm{ad}}}

\def\Sing{\mathop{\mathrm{Sing}}}
\def\i{\mathrm{i}}
\def\codim{\mathop{\mathrm{codim}}}

\def\d{\partial}
\def\la{\lambda}
\def\G{\mathit{\Gamma}}
\def\E{\mathcal{E}}
\def\R{\mathbb{R}}
\def\C{\mathbb{C}}
\def\P{\mathbb{P}}
\def\k{{\frak k}}
\def\g{{\frak g}}
\def\D{{\frak g}^*}
\def\O{\mathcal{O}}
\def\e{\eta}
\def\K{\mathcal{K}}
\def\S{\mathcal{S}}

\renewcommand{\ge}{\geqslant}
\renewcommand{\phi}{\varphi}

\begin{document}
\begin{frontmatter}
\title{Projections of Jordan bi-Poisson structures that are
Kronecker, diagonal actions, and the classical Gaudin
systems\thanksref{msc}}
\author{Andriy
Panasyuk}

\address{Division of Mathematical Methods in Physics,\\ University
of Warsaw,\\ Ho\.{z}a St.~74, 00-682
Warsaw, Poland,\\ e-mail:  panas@fuw.edu.pl\\ and \\
Pidstrygach Institute of Applied Problems of \\
Mathematics
and Mechanics,\\ Naukova Str. 3b,\\ 79601 Lviv, Ukraine}

\begin{abstract}
We propose a method of constructing completely integrable systems
based on reduction of bihamiltonian structures. More precisely, we
give an easily checkable necessary and sufficient conditions for
the micro-kroneckerity of the reduction (performed with respect to
a special type  action of a Lie group) of micro-Jordan
bihamiltonian structures whose Nijenhuis tensor has constant
eigenvalues. The method is applied to the diagonal action of a Lie
group $G$ on a direct product of $N$ coadjoint orbits
$\O=O_1\times\cdots\times O_N\subset\D\times\cdots\times\D$
endowed with a bihamiltonian structure whose first generator is
the standard symplectic form on $\O$. As a result we get the so
called classical Gaudin system on $\O$. The method works for a
wide class of Lie algebras including the semisimple ones and for a
large class of orbits including the generic ones and the
semisimple ones.
\end{abstract}

\end{frontmatter}
\date{}

\section{Introduction}
In this paper we propose a method of constructing completely
integrable systems, based on reduction of bihamiltonian
structures. The method is illustrated by producing a class of
systems on products of coadjoint orbits, which include the
so-called classical Gaudin systems. Now we will briefly explain
the method and discuss the Gaudin systems.

According to the last decade investigations of bihamiltonian
structures \cite{gz1,gz2,gz3,z1}, i.e. pairs of compatible Poisson
bivectors which will be called bi-Poisson structures in this
paper, there are two main classes of them: micro-Jordan and
micro-Kronecker (we shall omit the prefix micro- in this
Introduction for shortness). Jordan structures can be
characterized by the property that almost every bivector in the
corresponding pencil is nondegenerate, that is, the structure can
be generated by the inverses of two symplectic forms
$(\omega_1)^{-1},(\omega_2)^{-1}$. Kronecker pencils consist of
degenerate bivectors and are distinguished by the condition of the
constancy of rank (see Section \ref{s20}). Both classes of
bi-Poisson structures play important role in completely integrable
systems. Given a Jordan structure, one constructs an involutive
family of functions by means of the corresponding Nijenhuis
operator $N=\omega_2^{-1}\circ\omega_1$ (the eigenvalues of $N$
are in involution and in various examples there are enough
functionally independent ones); in Kronecker case functions in
involution appear as Casimir functions of the Poisson bivectors of
the pencil and form a complete set.

Among Jordan bi-Poisson structures there are ones which are
trivial from the point of view of complete integrability:
structures with the constant eigenvalues of $N$. We call them dull
(after I.S.Zakharevich). It is amazing that using the simultaneous
Poisson reduction of the symplectic forms $\omega_1, \omega_2$
related to such a structure one can produce a Kronecker bi-Poisson
structure which is far from being "dull" since it gives a complete
involutive family of functions, Casimirs of the Kronecker pencil.
Due to this remark the following question seems to be important:
when a simultaneous Poisson reduction of two symplectic forms
generating a dull Jordan bi-Poisson structure gives a Kronecker
one? We give necessary and sufficient conditions for such a
reduction to be Kronecker in the situation which roughly can be
described as follows: a Lie group $G$ acts freely on a manifold
$M$ with a dull Jordan bi-Poisson structure; this action is
hamiltonian with respect to all bivectors of the pencil; the
induced actions on the spaces of symplectic leaves of the
exceptional (i.e. of nonmaximal rank) bivectors are transitive
(see Theorem \ref{main1}).

Now, assume we are in such a situation and the necessary and
sufficient conditions mentioned are satisfied. Then we are able to
produce two complete involutive families of functions on $M/G$ and
$M$ respectively. First of them, $\mathcal{F}$, is generated by
all Casimir functions of bivectors from the constructed Kronecker
pencil on $M/G$. It is involutive and complete with respect to any
Poisson bivector of the pencil. The second one,
$\mathcal{G}^{t_0}$, is related to any nondegenerate bivector
$\e^{t_0}$ from the initial dull Jordan pencil
$\{\e^t=(\omega_1)^{-1}+t(\omega_2)^{-1}\mid t\in
\P^1=\R^1\cup\infty\}$. Denoting by $p$ the canonical projection
$M\rightarrow M/G$, we define $\mathcal{G}^{t_0}$ as the family
$p^*\mathcal{F}$ completed by $\mu_{t_0}^*\mathcal{F}'$, where
$\mu_{t_0}$ is the corresponding moment map $M\rightarrow\g^*$ and
$\mathcal{F}'$ is a complete involutive set of functions on $\g^*$
(endowed with the canonical linear Poisson bivector
$\e_{\mathrm{can}}$). The family $\mathcal{G}^{t_0}$ is involutive
and complete with respect to $\e^{t_0}$.

Note, that due to the standard properties of the dual pairs of
Poisson structures (cf. Section \ref{s10}) the family
$\mathcal{G}^{t_0}$ can be also generated by
$\mu_{t_0}^*\mathcal{F}'$ and by
$\{\mu_t^*Z^{\e_{\mathrm{can}}}\mid t\in \P^1\}$ (instead of
$p^*\mathcal{F}$), where $Z^{\e_{\mathrm{can}}}$ is the set of
Casimirs of $\e_{\mathrm{can}}$, i.e. invariants of coadjoint
representation.

Next, we apply the method described above to the following data:
$M=O_1\times\cdots\times O_N\subset (\g^*)^{\times N}$ is a
coadjoint orbit of the Cartesian product $G^{\times N}$ of $N$
copies of a Lie group $G$, $G$ acting on $M$ diagonally;
$\omega_1=\omega_{(1)}+\cdots+\omega_{(N)}$ is the standard
symplectic form on $M$, $\omega_{(i)}$ being the standard
symplectic form on $O_i$; $\omega_2$ is defined as
$(1/a_1)\omega_{(1)}+\cdots+(1/a_N)\omega_{(N)}$, where
$a_1,\ldots,a_N$ are any different real numbers. Under some
conditions on the orbits $O_1,\ldots, O_N$ (see Theorem
\ref{main2}) Theorem \ref{main1} can be applied and we get a
Kronecker bi-Poisson structure on the regular part of the variety
$M/G$ and the corresponding complete involutive sets of functions
$\mathcal{F}$ and $\mathcal{G}^{t_0}$. By the remark above this
last can be generated by $\mu_{t_0}^*\mathcal{F}'$ and
$\{\mu_t^*Z^{\e_{\mathrm{can}}}\mid t\in \P^1\}$, where
$\mu_t(x_1,\ldots,x_N)= (1/(t+a_1))x_1+\cdots+(1/(t+a_N))x_N$, as
calculations show. So, we recognize in $\mathcal{G}^{t_0}$ the so
called classical Gaudin integrable system.

The quantum version of this system, which uses the Lie algebra
$\mathrm{su(2)}$ and describes some type of interaction of
particles with spin, was introduced by M.Gaudin \cite{g1,g2,g3}.
Later E.Sklyanin studied separability of classical and quantum
systems in case of $\g=\mathrm{sl(n)}$ and with additional term in
$\mu_t$, a constant regular matrix \cite{s1,s2,s3}. The
integrability of such systems was also discussed in \cite{rst}
from the point of view of r-matrix formalism.

Summarizing, the following items of this paper seem to be new: 1)
the method of constructing completely integrable systems based on
the reduction of dull Jordan bi-Poisson structures; 2) application
of this method to the Gaudin type systems; 3) proof of the
complete integrability of such systems for a wide class of Lie
algebras including all semisimple ones and for a wide class of
coadjoint orbits including all generic ones and all semisimple
ones (see Section \ref{s40}, in particular Remark \ref{40.90}).

The paper is organized as follows. Section \ref{s10} is
preparatory: we introduce notations and recall standard
definitions and facts related to Poisson structures, their dual
pairs and hamiltonian actions. Proposition \ref{10.140} is new
(but easy). Proposition \ref{10.150} is doubtlessly known but the
author was not able to find its explicit formulation in the
literature.

Similarly, Section 2 serves for introducing the notations and main
data on bi-Poisson structures. The material of this section is
more or less standard.

In Section \ref{s30} we give  the first main result of this paper:
necessary and sufficient conditions for the kroneckerity of the
reduction of a dull Jordan bi-Poisson structure with respect to a
specific bi-hamiltonian action of a Lie group (Theorem
\ref{main1}). In Corollary \ref{30.cor} we apply this result for
constructing a completely integrable system on the initial
manifold (the above mentioned family $\mathcal{G}^{t_0}$). We also
illustrate the method by an example of a diagonal action of
$\mathrm{SL(2)}$ on $\R^{2N}$ endowed with a dull Jordan
bi-Poisson structure (see \ref{30.exa}).

In Section \ref{s40} we develop this example and construct a dull
Jordan  bi-Poisson structure on a cartesian product of $N$
coadjoint orbits of a Lie group $G$, whose reduction with respect
to the diagonal action of $G$ is Kronecker. Theorem \ref{main2}
which establishes this kroneckerity using Theorem \ref{main1} is
the second main result of the paper. Also, we calculate the
corresponding families of the moment maps $\mu_t$ and complete
involutive families of functions $\mathcal{F}$ and $\mathcal{G}^t$
(Corollaries \ref{40.60}--\ref{40.80}). We conclude the paper by
the discussion on the range of applicability of the method
(Subsections \ref{dissc}--\ref{40.100}).

\section{Projections of Poisson structures, dual pairs and
complete involutive sets of functions}
\label{s10}
\abz\label{10.05}
\begin{conv}\rm
All objects in this paper are real-analytic or complex analytic,
$M$ stands for a connected manifold, $\E(U)$ for a space of
respectively real-valued analytic or holomorphic functions on an
open set $U\subset M$. We shall write $\mathbb{K}$ for $\R$ or
$\C$ depending on the category. The terms reduction and projection
related to the Poisson structures are synonyms in this paper.
\end{conv}
\abz\label{10.10}
\begin{defi}
Let $M$ be a manifold, $\eta\in\G(\bigwedge^2TM)$ be a bivector
field (from now on we shall skip the last word). We consider $\e$
as a homomorphism \[ \e:T^*M\longrightarrow TM
\]
obtained by the contraction on the first index and define
the (generalized) distribution of characteristic subspaces
$\chi^\e\subset TM$ by
\[\chi^\e_x=\im \e_x,\,\,\, x\in M.
\]
Set $\rk\e_x=\dim\im\e_x, \rk\e=\max_{x\in M}\rk\e_x$ and
$R^\e=\{x\in M\mid\rk\e_x=\rk\e\}$.
We say
that $\e$ is nondegenerate if it is isomorfizm or,
equivalently, $\chi^\e=TM$.
\end{defi}
Clearly, $R^\e\subset M$ is an open dense set.
\abz\label{10.20}
\begin{defi}
Let $M$ be a manifold, $\K$ be a foliation on $M$ such that
the factor space $M'=M/\K$ is a manifold, and let
$p:M\rightarrow M'$ be the canonical projection. We say
that a bivector $\e\in\G(\bigwedge^2TM)$ is projectable via
$p$ if there exists a bivector $\e'\in\G(\bigwedge^2TM')$
(called the projection of $\e$) such that
\[
\e'_{x'}=p_*\e_x
\]
for any $x'\in M', x\in p^{-1}(x')$.
\end{defi}
\abz\label{10.40}
\begin{defi}
A bivector $\e\in\G(\bigwedge^2TM)$ is called Poisson if
the operation
\[
\{f,g\}^\e=\e(f)g,\ \ f,g\in \E(M),
\]
where we put $\e(f)=\e(df)$, satisfies the Jacobi identity.
The operation $\{\,,\}^\e$ is called the Poisson
bracket, the vector fields $\e(f)$ are called hamiltonian.
\end{defi}
\abz\label{10.45}
\begin{propo}(e.g. \cite{lichn})
A bivector $\e$ is Poisson iff $[\e,\e]=0$, where $[\,,]$
is the Schouten bracket on mutivector fields.
\end{propo}
\abz\label{10.50}
\begin{theo} (\cite{kir})
If $\e$ is a Poisson bivector its generalized distribution of
characteristic subspaces $\chi^\e$ is completely integrable, i.e.
there exists a generalized foliation $\mathcal{S}$ on $M$ such
that $T_x\mathcal{S}=\chi^\e_x$ for any $x\in M$. The restriction
$\e|_S$ of $\e$ to any leaf $S$ of $\mathcal{S}$ is a correctly
defined nondegenerate Poisson bivector.
\end{theo}
\abz\label{10.60}
\begin{defi}
The leaves of the generalized foliation $\mathcal{S}$ are called
symplectic leaves of $\e$.
\end{defi}
The definition is motivated by the fact that the inverse to a
nondegenerate Poisson bivector 2-form is symplectic.
\abz\label{10.70}
\begin{defi}
A function $f\in \E(U)$ over an open set $U\subset
M$ is called a Casimir function for $\e$ if $\e(f)\equiv
0$. The set of all Casimir functions for $\e$ over $U$ will
be denoted by $Z^\e(U)$.
\end{defi}
Geometrically speaking the Casimir functions are those constant on
the symplectic leaves of maximal dimension.
\abz\label{10.75}
\begin{defi}
A set $Z\subset Z^\e(U)$ of Casimir functions over $U\subset M$ is
called complete if there exist $f_1,\ldots,f_k\in Z$, where
$k=\crk\e:=\dim M-\rk\e$, such that their differentials are
independent on $U\cap R^\e$ and for any $F\in \E(\mathbb{K}^k)$
the composition $F(f_1,\ldots,f_k)$ belongs to $Z$.
\end{defi}
In other words $Z$ is complete iff the common level sets of
functions from $Z$ coincide with the symplectic foliation
on $U\cap R^\e$.

It is clear
that $Z^\e(U)$ is complete for sufficiently small $U$.
\abz\label{10.80}
\begin{defi}
A set $I\subset  \E(U)$ of functions over $U\subset M$ is called
complete involutive for $\e$ if: 1) $\{f,g\}^\e=0\ \ \ \forall f,g
\in I$; 2) there exist $f_1,\ldots,f_s\in I$, where $s=\dim
M-(1/2)\rk\e$, such that their differentials are independent on
$U\cap R^\e$ and for any $F\in \E(\mathbb{K}^s)$ the composition
$F(f_1,\ldots,f_s)$ belongs to $I$.
\end{defi}
If $I$ is a complete involutive set over $U$, then $I\supset
Z^\e(U)$ and the last set is complete in the sense of \ref{10.75}.
Any such set $I$ is a set of functions constant on a foliation of
$U\cap R^\e$ of dimension $(1/2)\rk\e$ which is lagrangian in  any
symplectic leaf (of maximal dimension).
\abz\label{10.84}
\begin{defi}
A map $\mu:(M,\e)\rightarrow (M,',\e')$ between two Poisson
manifolds is called Poisson if for any $f,g\in\E(M')$
\[
\mu^*\{f,g\}^{\e'}=\{\mu^*f,\mu^*g\}^\e
\]
or, equivalently, $\mu_*\e_x=\e'_{\mu(x)}$ for any $x\in M$.
\end{defi}
\abz\label{10.85}
\begin{propo}(\cite{weinst}, Lemma 1.2)
If $\mu:(M,\e)\rightarrow (M_1,\e_1)$ is a Poisson map, then the
trajectory of any hamiltonian field $\e_1(f)$ is the projection
via  $\mu$ of the trajectory of hamiltonian field $\e(\mu^*f)$.
\end{propo}

\abz\label{10.86}
\begin{propo}
Let $p:M\rightarrow M'$ be as in \ref{10.20} and assume that $\e$
is a Poisson bivector on $M$. Then the bivector $\e$ is
projectable via $p$ iff for any open set $U\subset M'$ the
subspace $p^*\E(U)=\{p^*f\mid f\in \E(U)\}\subset \E(P^{-1}(U))$
is a Lie subalgebra with respect to $\{\,,\}^{\e|U}$. If $\e$ is
projectable and $\e'$ is the projection, then $\e'$ is a Poisson
bivector and  $p:(M,\e)\rightarrow (M',\e')$ is a Poisson map.
\end{propo}

\begin{pf} Let $(U,\{\phi_j\})$ be a coordinate map on $M'$. Since
$p^*\E(U)$ is a subalgebra, $\{p^*\phi_i,p^*\phi_j\}^\e=p^*c^{ij}$
for some function $c^{ij}\in\E(U)$. It is easily seen that
$c^{ij}$ transfoms tensorially under coordinate changes, i.e.
represents some bivector $\e'$ on $M'$. The remaining part of the
proof is almost immediate consequence of the definitions.
\qed\end{pf}

Here is another criterion of projectability.
\abz\label{10.90}
\begin{theo}
(Liebermann-Weinstein criterion of projectability,
\cite{lieb,weinst}) Let $p:M\rightarrow M'$ and $\K$ be as in
\ref{10.20} and let $\e$ be a nondegenerate Poisson bivector on
$M$. Write $N\K\subset T^*M$ for the conormal bundle to the
foliation $\K$. Then $\e$ is projectable via $p$ iff the
distribution $\e(N\K)\subset TM$, which is the skew-orthogonal
complement to the distribution $T\K$, is completely integrable.
\end{theo}
\abz\label{10.100}
\begin{coro}(\cite{weinst})
Let $p:M\rightarrow M'$, $\K$, and $\e$ be as in the assumption of
Theorem \ref{10.90}. Assume that $\e$ is projectable and that the
foliation $\K''$ tangent to the distribution $\e(N\K)$ is such
that the factor space $M''=M/\K''$ is a manifold. Then $\e$ is
also projectable to $M''$ via the canonical projection
$p'':M\rightarrow M''$.
\end{coro}

\begin{pf*}{Proof} follows from the fact that in the nondegenerate case
$\e(N[\e(N\K)])=T\K$, i.e. the distributions $T\K$ and $\e(N\K)$
are the skew-orthogonal complements of each other. \qed\end{pf*}
\abz\label{10.110}
\begin{defi} (\cite{weinst})
Let $\e$ be a nondegenerate Poisson bivector on $M$ and let
$\K',\K''$ be foliations on $M$
such that $\e(N\K')=T\K''$ and the factor spaces
$M'=M/\K', M''=M/\K''$ are manifolds. The pair
$(\e',\e'')$, where $\e'=p'_*\e, \e''=p''_*\e$ are the
projections of $\e$ via the canonical projections
$p':M\rightarrow M'$ and $p'':M\rightarrow M''$
respectively, is called a dual pair of Poisson bivectors.
\end{defi}

The situation can be expressed by the following diagram
\begin{center}
\setlength{\unitlength}{1mm}
\begin{picture}(90,33)(0,0)
\put(45,30){\makebox(0,0){$(M,\e)$}}
\put(22,5){\makebox(0,0){$(M'=M/\K',\e'=p'_*\e)$}}
\put(70,5){\makebox(0,0){$(M''=M/\K'',\e''=p''_*\e),$}}
\put(24,18){\makebox(0,0){$(p',p'_*)$}}
\put(68,18){\makebox(0,0){$(p'',p''_*)$}}
\put(42,26){\vector(-1,-1){17}}
\put(48,26){\vector(1,-1){17}}
\end{picture}
\end{center}
where $\e(N\K')=T\K''$.
\abz\label{10.115}
\begin{exa}\rm
Let $G$ be a connected Lie group with the Lie algebra $\g$.
Assume it is acting on a Poisson manifold $(M,\e)$, in
particular a Lie algebra homomorphism $\rho:\g\rightarrow
\G TM$ is given (the space of vector fields is endowed with
the commutator Lie bracket).

The action is called hamiltonian if
there exists a Lie algebra homomorphism $\psi:\g\rightarrow
\E (M)$ such that the following diagram is commutative
\[
\begin{array}{ccc}
\g & \stackrel{\psi}{\longrightarrow}&
\E(M) \\ \parallel & & \downarrow  \e(\cdot) \\ \g
&\stackrel{\rho}{\longrightarrow} & \G TM,
\end{array}
\]
where $\e(\cdot)$ is the Lie algebra homomorphism of taking the
hamiltonian vector field (see \ref{10.40}). The map $\mu: x\mapsto
\varphi_x:M\rightarrow \D$ defined by $\varphi_x(v)=\psi(v)(x),
v\in\g$ is called the moment map.

Assume that G acts on $M$ by Poisson maps and that $M/G$ is a
manifold. Then by Proposition \ref{10.86} the bivector $\e$ is
projectable via the canonical projection $M\rightarrow M/G$. If
moreover, $\e$ is nondegenerate and the action is hamiltonian its
orbits are skew-orthogonal to the fibers of the moment map. This
last is a Poisson map from $(M,\e)$ to $(\D,\e_{\mathrm{can}})$,
where $\e_{\mathrm{can}}$ is the canonical linear Poisson bivector
on $\D$. In case of a locally free action $\mu$ is a submersion
(see the proof of Corollary \ref{10.130}) and we get a dual pair
$(\e',\e_{\mathrm{can}})$.
\end{exa}

We complete the section by a series of propositions that
will be crucial for the subsequent part of the paper.
Proposition \ref{10.118} and Corollaries
\ref{10.120}, \ref{10.130} are classical.
\abz\label{10.118}
\begin{propo}
We retain the notations of \ref{10.110}. Assume that $(\e',\e'')$
is a dual pair of Poisson bivectors. Then
\begin{enumerate}
\item[(a)] the distribution
$D_x=T_x\K'+T_x\K''\subset T_xM, x\in M$, is of constant dimension
on an open dense set $R\subset M$ and is completely integrable on
$R$;
\item[(b)] the foliation $\mathcal{D}$ tangent to $D$ on $R$ is the pull-back
of the foliations of symplectic leaves $\mathcal{S}',
\mathcal{S}''$ of maximal dimension of the bivectors $\e',\e''$
respectively:
\[
(p')^*\mathcal{S}=\mathcal{D}=(p'')^*\mathcal{S}'';
\]
\item[(c)] $\crk\e'=\crk\e''$.
\end{enumerate}
\end{propo}
\begin{pf} The constancy of dimension on an open dense set follows from
analyticity of all objects. Item (b) is a consequence of
Proposition \ref{10.85} and of skew-orthogonality of the
foliations $\K',\K''$. Also, (b) implies integrability of $D$ and
(c). \qed\end{pf}

\abz\label{10.120}
\begin{coro}
We retain the assumptions of \ref{10.118}. Let $U'\subset M',
U''\subset M''$ be open sets such that the  sets of Casimir
functions $Z':=Z^{\e'}(U'), Z'':=Z^{\e''}(U'')$ are complete (see
\ref{10.75}) and $U:=(p')^{-1}(U')\cap (p'')^{-1}(U'')\not
=\emptyset$. Put $((p')^*Z')|_U=\{((p')^*f)|_U \mid\ f\in Z'\}$
and $((p'')^*Z'')|_U=\{((p'')^*g)|_U \mid g\in Z''\}$. Then \sabz
\begin{equation}
((p')^*Z')|_U=((p'')^*Z'')|_U.       \label{f10.110}
\end{equation}
\end{coro}
\abz\label{10.130}
\begin{coro}
Assume that a Lie group $G$ is acting in the hamiltonian way on a
symplectic Poisson manifold $(M,\e)$ (see Example \ref{10.115}).
Assume moreover, that this action is locally free (the stabilizer
of any point is at most discrete) and that $M/G$ is a manifold.
Then for any $x'\in M/G$ we have $\crk\e_{x'}'=\rk G$, where $\e'$
is the projection of $\e$ via the canonical map $M\rightarrow
M/G$.
\end{coro}

\begin{pf} It is well known that the image of the differential at a point
$x\in M$ of the moment map $\mu:M\rightarrow \D$ coincides with
the annihilator in $\D$ of $\g_x\subset \g$, where $\g_x$ is the
Lie algebra of the stabilizer of $x$ (see \cite{gs}, Lemma 2.1).
Thus in our situation when the stabilizer is discrete $\mu$ is a
submersion. By \ref{10.118} corank of $\e'$ coincides with the one
of $\e_{\mathrm{can}}$, i.e. with rank of the Lie group $G$.
\qed\end{pf}

Here is a generalization of this result to the case of
degenerate Poisson bivector $\e$.
\abz\label{10.140}
\begin{propo}
Let $\e$ be a regular Poisson bivector (i.e. $\rk\e_x=const$) on
$M$ and let a Lie group $G$ act locally freely on $M$ in such a
way that $M/G$ is a manifold. Given a symplectic leaf $S\subset
M$, write $G_S\subset  G$ for its stabilizer, i.e. for a subgroup
defined by $G_SS\subset S$. Fix $S$ and assume that
\begin{enumerate}
\item[(1)] $G$ acts by Poisson maps, i.e. the action preserves $\e$;
\item[(2)] the action induces
a transitive action on the space of symplectic leaves;
\item[(3)] the induced action of $G_S$ on $(S,\e|_S)$ is
hamiltonian.
\end{enumerate}
Then
\begin{enumerate}
\item[(a)] if $\widehat{S}$ any symlectic leaf, the
stabilizers $G_S, G_{\widehat{S}}$ are conjugate;
\item[(b)] the induced action of $G_{\widehat{S}}$ on
$\widehat{S}$ is hamiltonian;
\item[(c)] $\e$ is projectable via the
canonical map $M\rightarrow M/G$ and $\crk\e_{x'}'=\rk G_S$
for any $x'\in M/G$, where $\e'$ is the projection.
\end{enumerate}
\end{propo}

\begin{pf} Since any two points on any  symplectic leaf $S$ of $\e$ can
be connected by a finite number of hamiltonian trajectories and
since the action preserves $\e$, it follows from Proposition
\ref{10.85} that the image $gS, g\in G$, is again a symplectic
leaf. Now, assumption (2) implies that for any $S,\widehat{S}$
there exists $a\in G$ such that $aS=\widehat{S}$, hence
$G_S=\{g\in G\mid gS= S\}=\{g\in G\mid
ga^{-1}\widehat{S}=a^{-1}\widehat{S}\}= \{g\in G\mid
aga^{-1}\widehat{S}=\widehat{S}\}=a^{-1}G_{\widehat{S}}a$.

To prove (b) let us consider the induced action
$\rho_{\widehat{S}}:\g_{\widehat{S}}\rightarrow \G T\widehat{S}$
of the Lie algebra of the stabilizer $G_{\widehat{S}}$ on
$\widehat{S}$. Its hamiltonicity follows from the following
commutative diagram:
\[
\begin{array}{lclcl}
\g_{\widehat{S}} &
\stackrel{\rho_{\widehat{S}}}{\longrightarrow}&
\G T\widehat{S} & = & \G T\widehat{S} \\
\downarrow \Ad a &  & \uparrow L_{a*}  & & \parallel\\
\g_S & \stackrel{\rho_S}{\longrightarrow} & \G TS & &
\G T\widehat{S}\\
\parallel & & \uparrow  \e(\cdot) &  & \uparrow
\e(\cdot)\\
\g_S & \stackrel{\psi_S}{\longrightarrow} & \E(S) & &
\E(\widehat{S})\\
\parallel & & \uparrow L_a^* & & \parallel\\
\g_S & \stackrel{L_{a^{-1}}^*\circ\psi_S}{\longrightarrow} &
\E(\widehat{S}) & = & \E(\widehat{S}),
\end{array} \]
where all maps are Lie algebra homomorphisms, $\psi_S$ is one
existing by assumption (3), $L_a$ denotes the left multiplication
by $a$.

Projectability of $\e$ follows from (1) and from \ref{10.86}.
Condition (2) guarantees that the projection $\e'$ of $\e$ via the
map $M\rightarrow M/G$ coincides with the projection $(\e|_S)'$ of
the resricted Poisson bivector $\e|_S$ via the map $S\rightarrow
S/G_S=M/G$.  Taking into account assumption (3) we can apply
Corollary \ref{10.130} to the action of $G_S$ on $(S,\e|_S)$. This
proves (c). \qed\end{pf}

\abz\label{10.150}
\begin{propo}
We retain the notations of \ref{10.110}. Let $(\e',\e'')$ be a
dual pair of Poisson bivectors, let $U'\subset M', U''\subset M''$
be open sets such that $U:=(p')^{-1}(U')\cap (p'')^{-1}(U'')\not
=\emptyset$ and let $I'\subset \E(U'), I''\subset \E(U'')$ be
complete involutive sets of functions for $\e', \e''$
respectively. Put $((p')^*I')|_U=\{((p')^*f)|_U \mid f\in I'\}$
and $((p'')^*I'')|_U=\{((p'')^*g)|_U \mid g\in I''\}$. Then the
space $I:=((p')^*I')|_U+((p'')^*I'')|_U$ is a complete involutive
set of functions for $\e$.
\end{propo}

\begin{pf} We first notice that since $\K'$ and $\K''$
are skew-orthogonal, $\{(p')^*f,\lbr(p'')^*g\}^\e=0$ for any $f\in
\E(U'), g\in \E(U'')$. Together with the Poisson property for $p'$
and $p''$ this shows that $I$ is an involutive set of functions
with respect to $\e$. Now we only need to calculate its
"functional dimension".

Let us choose a "functional basis" $\{f_1,\ldots,f_{s'}\}$
of $I'$ such that $f_1,\ldots,f_{r'}\in Z^{\e'}(U')$ and
any "functional basis" $\{g_1,\ldots,g_{s''}\}$ of $I''$.
Then the functions
$(p')^*f_{r'+1},\ldots,(p')^*f_{s'},(p'')^*g_1,\ldots,(p'')^*g_{s''}$
are functionally independent on an open dense subset of $U$
since
\[
\{(p')^*f|_U\mid f\in \E(U')\}\cap
\{(p')^*g|_U\mid g\in \E(U'')\}=Z,
\]
where $Z$ denotes the set (\ref{f10.110}). Now, one has
\begin{eqnarray*}
s'-r'&=&\frac{1}{2}\rk\e'=\frac{1}{2}(\dim\K''-\dim\K''\cap\K'),\\
s''&=&\frac{1}{2}\rk\e''+\crk\e''=\frac{1}{2}(\dim\K'-\dim\K''\cap\K')+
\dim\K''\cap\K',
\end{eqnarray*}
and, finally
\[
s'-r'+s''=\frac{1}{2}(\dim\K''+\dim\K')=\frac{1}{2}\dim M.
\]
\qed\end{pf}

\section{Preliminaries on bi-Poisson structures}
\label{s20}
\abz\label{20.10}
\begin{defi}
A pair $(\e_1,\e_2)$ of linearly independent bivectors on a
manifold $M$ is called Poisson if $\e^t:=t_1\e_1+t_2\e_2$ is a
Poisson bivector for any $t=(t_1,t_2)\in \mathbb{K}^2$; the whole
family of Poisson bivectors $\{\e^t\}_{t\in\mathbb{K}^2}$ is
called a bi-Poisson structure.  By definition the family
consisting of zero bivector is a bi-Poisson structure called
trivial.
\end{defi} A bi-Poisson structure $\{\e^t\}$ (we shall
often skip the parameter space in the notations) can be
viewed as a two-dimensional vector space of Poisson
bivectors, the Poisson pair $(\e_1,\e_2)$ as a basis in
this space. Of course, the basis can be changed.
\abz\label{20.20}
\begin{defi}
A bi-Poisson structure $\{\e^t\}$ is called Jordan at a
point $x\in M$ if $\rk\e_x^t=\dim M$ for some $t$.
A bi-Poisson structure is called micro-Jordan if it is
Jordan at any point of some open dense subset in $M$.
\end{defi}
The terminology is due to I.Gelfand and I.Zakharevich
\cite{gz1,z1} who reduced the analysis of a bi-Poisson structure
at a point to the study of a pencil of operators and applied the
classical classificational results. These last say that any pencil
is built of the irreducible ones, the so-called Jordan and
Kronecker blocks. The above definition corresponds to the case
when only the Jordan blocks are present.

The theory of pencils of operators is well understood over
the field of complex numbers. We shall also need some
notions related to the complexification matters.
\abz\label{20.30}
\begin{nota}
\rm If $M$ is a real manifold (recall that objects are real
analytic) we denote by $\widetilde{M}$ some
complexification of $M$, i.e. a complex manifold
$\widetilde{M}$ such that $M$ is embedded in $\widetilde{M}$ as a
totally real submanifold. The complex structure near $M$ is
defined uniquely up to a biholomorphic map preserving $M$
 (see \cite{bw}), that's why we use the same notation
$\widetilde{M}$ for possibly different complexifications.
Given any tensor $\e$ on $M$, we write $\widetilde{\e}$ for
its complexification, which is a holomorphic
tensor defined on $\widetilde{M}$ (the last should be
shrinked if needed).

For any real bi-Poisson structure $\{\e^t=t_1\e_1+t_2\e_2\}$ on
$M$ we denote by $\widetilde{\e}^t$ its complexification, i.e. the
holomorphic bi-Poisson structure
$\{\widetilde{\e}^t=t_1\widetilde{\e}_1+t_2\widetilde{\e}_2\mid
t=(t_1,t_2)\in\C^2\}$ on $\widetilde{M}$.

If $M$ and $\{\e^t\}$ are a priori holomorphic we put
$\widetilde{M}=M$, $\{\widetilde{\e}^t\}=\{\e^t\}$ e.t.c., hence
tilde for holomorphic objects  will denote themselves (not the
complexification of the underlying real-analytic objects).
\end{nota}
\abz\label{20.40}
\begin{defi}
Let $\{\e^t\}$ be a micro-Jordan bi-Poisson structure on $M$. Put
$E(x)=\{t\in\C^2\mid\rk_\C\widetilde{\e}_x^t<\dim_\C
\widetilde{M}\}\subset\C^2, x\in M$. This set is called
exceptional for $\{\e^t\}$ at $x$. If $E=E(x)$ does not depend on
$x$ the structure $\{\e^t\}$ is called dull.
\end{defi}
This terminology is due
to I.Zakharevich and is motivated by the fact that the
constancy of $E(x)$ implies the constancy of the
eigenvalues for the recursion operator
$\e_1^{-1}\circ\e_2$, i.e. the situation is far from being
of interest in the theory of integrable systems in which
these eigenvalues appear as the first integrals.

It is clear that $E(x)$ consists of a finite number of
1-dimensional subspaces in $\C^2$.
\abz\label{20.50}
\begin{defi}
Let $\{\e^t\}$ be a bi-Poisson structure on $M$. It is
called Kronecker at a point $x\in M$ if $\rk_\C\widetilde{\e}_x^t$
is constant with respect to $t\in\C^2\setminus\{0\}$.
We say that $\{\e^t\}$ is micro-Kronecker if it is
Kronecker at any point of some open dense set in $M$.
In particular the trivial bi-Poisson structure is
micro-Kronecker.
\end{defi}
Again this terminology is due to I.Zakharevich and is
motivated by the fact that under the above rank assumptions
the corresponding pencil of operators (see the discussion
in \ref{20.20}) contains only the Kronecker blocks.

\abz\label{30.10}
\begin{defi}
Let $p:M\rightarrow M'$ be as in \ref{10.20} and let
$\{\e^t=t_1\e_1+t_2\e_2\}$ be a bi-Poisson structure on $M$.  We
say that it is projectable via $p$ if so is the bivector $\e^t$
for any $t$. The family $\{(\e^t)'=t_1\e'_1+t_2\e'_2\}$ consisting
of the projections of $\e^t$, which is a bi-Poisson structure on
$M'$ under the condition that the bivectors $\e'_1,\e'_2$ are
linearly independent or trivial (see Propsition \ref{10.86}), is
called the projection of $\{\e^t\}$.
\end{defi}

Now we are able to formulate the main question of this
paper: when the projection of a (projectable) dull
micro-Jordan bi-Poisson structure is micro-Kronecker?
We shall
answer it in the next section for some particular cases of
locally free bi-Poisson actions.  Now we want to present a
result which shows why the micro-Kronecker structures are
interesting and which will be effectively used later.
\abz\label{20.60}
\begin{propo}
Let $\{\e^t\}$ be a micro-Kronecker bi-Poisson structure on $M$.
Assume that an open set $U\subset M$ is such that the set
$Z^{\e^t}(U)$ of Casimir functions for $\e^t$ over $U$ is complete
(see Definition \ref{10.75}) for any $t\not=0$. Then the set
\[
Z^{\{\e^t\}}(U):=\sum_{t\not=0}Z^{\e^t}(U)
\]
is a complete involutive set of functions for any $\e^t\not =0$
(see Definition \ref{10.80}). (Here and subsequently in similar
situations we understand the sum as the algebraic sum of linear
(sub)spaces of functions in the linear space of all functions. In
other words this sum coincides twith the linear span $\langle
Z^{\e^t}(U)\mid t\not =0\rangle$. Of course, it is enough to sum
over a sufficiently large finite set of indices $t$.)
\end{propo}
We shall call the functions from  $Z^{\{\e^t\}}$ the first
integrals of the bi-Poisson structure $\{\e^t\}$. The
reader is referred to a celebrated paper of A.Bolsinov
\cite{bols} for the proof of completeness.  Although the
involutivity of this set was known and extensively used
since the 80-ies the author was not able to find its proof
and gave a version of it in \cite{p}.
\abz\label{20.70}
\begin{exa}\rm
(Method of the argument translation) Let $\g$ be a Lie algebra
with $\codim \Sing\g^*\ge 3$, where $\Sing\g^*\subset \D$ is the
algebraic set of all coadjoint orbits of nonmaximal dimension (in
particular $\g$ can be any semisimple). Let
$\e_1=\e_{\mathrm{can}}$ be a canonical linear Poisson bivector on
$\D$, and let $\e_2=\e_{\mathrm{can}}(a)$ be the Poisson bivector
obtained by "freezing" $\e_{\mathrm{can}}$ at a regular (i.e.
belonging to $\D\setminus\Sing\D=R^{\e_{\mathrm{can}}}$) element
$a$. It is well-known that $(\e_1,\e_2)$ is a Poisson pair and
that the corresponding bi-Poisson structure
$\{\e^t_{\mathrm{AT}}\}$ is micro-Kronecker (see
\cite{bols,p,z1}). The set of first integrals
$Z^{\{\e^t_{\mathrm{AT}}\}}$ is functionally generated by
$f_1(x+\la a),\ldots,f_k(x+\la a), \la \in\mathbb{K}$, where
$f_1,\ldots ,f_k$ are the invariants of the coadjoint action.
\end{exa}

\section{A locally free bi-Poisson action of
a Lie group on a dull micro-Jordan structure}
\label{s30}
\abz\label{a30.10}
\begin{assu}\rm
Let $G$ be a real Lie group. We shall assume that it possesses the
complexification, i.e. a complex Lie group $\widetilde{G}=G^\C$
containing $G$ as a real subroup such that its Lie algebra $\g^\C$
is the complexification of the Lie algebra $\g$ of $G$. In
particular, $G$ may be linear semisimple or compact.

Given a real dull micro-Jordan bi-Poisson structure
$\{\e^t=t_1\e_1+t_2\e_2\}$ on a manifold $M$,  we denote by
$\widetilde{M}$ a complexification of $M$ such that the bivectors
$\e_1,\e_2$ are exteded to holomorphic Poisson bivectors
$\widetilde{\e}_1,\widetilde{\e}_2$ (automatically forming a
Poisson pair on $\widetilde{M}$). We write $\{\widetilde{\e}^t\}$
for the holomorphic bi-Poisson structure
$\{t_1\widetilde{\e}_1+t_2\widetilde{\e}_2\}$, $e_1,\ldots,e_N$
for the vectors in $\C^2$ spanning the lines of the exceptional
set $E=\langle e_1\rangle\cup\cdots\cup\langle e_N\rangle$ (see
Definition \ref{20.40}), and
$\widetilde{\e}^{e_1},\ldots,\widetilde{\e}^{e_N}$ for the
corresponding exceptional bivectors.

We retain the convention that $\widetilde{(\cdot)}=(\cdot)$ for a
holomorphic object $(\cdot)$ (cf. \ref{20.30}).

\end{assu}
The central result of this paper is the following.
\abz\label{main1}
\begin{theo}
We retain the above assumptions and notations. Assume a Lie group
$G$ is acting locally freely on a manifold $M$ with a dull
micro-Jordan bi-Poisson structure $\{\e^t\}$, that this action is
extended to a locally free action of $\widetilde{G}$ on
$\widetilde{M}$ (in the complex case this extended action is the
initial one) and that $M/G, \widetilde{M}/\widetilde{G}$ are
manifolds.  For any $j=1,\ldots,N$ fix a symplectic leaf $S_j$ of
maximal dimension of the exceptional bivector
$\widetilde{\e}^{e_j}$ and let $\widetilde{G}_j$ denote its
stabilizer.

We make the following additional
assumptions on the $\widetilde{G}$-action on
$\widetilde{M}$:
\begin{enumerate}
\item[(1)] it is
bi-Poisson, i.e. preserves $\widetilde{\e}_1,
\widetilde{\e}_2$;
\item[(2)] it induces a transitive action on the space of
symplectic leaves of maximal dimension of any exceptional
bivector $\widetilde{\e}^{e_j}$;
\item[(3)] the induced action of $\widetilde{G}_j$ on
$(S_j,\widetilde{\e}^{e_j}|_{S_j})$ is hamiltonian;
\item[(4)] the action of $\widetilde{G}$ on
$(\widetilde{M},\widetilde{\e}^t), t\in\C^2\setminus E$, is
also hamiltonian.
\end{enumerate}
Then \begin{itemize} \item $\{\e^t\}$ is projectable via the
canonical map $p:M\rightarrow M/G$; \item the projection
$\{(\e^t)'\}$ is a bi-Poisson structure under the condition that
the bivectors $\e'_1,\e'_2$ are linearly independent or trivial;
\item $\{(\e^t)'\}$ is Kronecker at any point $x'\in M/G\setminus
p(R^{\widetilde{\e}^{e_1}}\cup\cdots\cup
R^{\widetilde{\e}^{e_N}})$ iff
\[ \rk\widetilde{G}=\rk\widetilde{G}_1=\cdots=\rk\widetilde{G}_N
\]
(recall that $R^\e$ stands for the regularity set of a bivector
$\e$, see \ref{10.10}).
\end{itemize}
\end{theo}

\begin{pf} It is clear that each $\e^t$ is projectable (since $G$ acts by
the Poisson maps with respect to $\e_1,\e_2$, see Proposition
\ref{10.86}) and that $\{(\e^t)'\}$ is a bi-Poisson structure
provided $\e_1',\e_2'$ are linearly independent or trivial.

By definition $\{(\e^t)'\}$ is Kronecker at $x'$ iff
$\crk(\widetilde{\e}^t)'_{x'}$ is constant with respect to $t\not
=0$. Now it remains to use Corollary \ref{10.130} to deduce that
$\crk(\widetilde{\e}^t)'_{x'}=\rk \widetilde{G}$ for $t\in
\C^2\setminus E$ and Proposition \ref{10.140} to get
$\crk(\widetilde{\e}^{e_j})'_{x'}=\rk\widetilde{G}_j,
j=1,\ldots,N$. \qed\end{pf}
\abz\label{30.cor}
\begin{coro}
In the situation of the above theorem let $\mu_t:M\rightarrow
\g^*, t\in\mathbb{K}^2\setminus E$, denote the moment map
corresponding to $\e^t$. Assume that $\{(\e^t)'\}$ is Kronecker.
Then
\begin{enumerate}
\item[(a)]
the pullback of the set of first integrals
$p^*\mathcal{F}:=p^*(Z^{\{(\e^t)'\}})$ (see \ref{20.60}) is equal
to
\[
p^*\mathcal{F}=\sum_{s\in \mathbb{K}^2\setminus
E}\mu_s^*(Z^{\e_{\mathrm{can}}}),
\]
where $\e_{\mathrm{can}}$ is the canonical linear Poisson bivector
on the dual space $\g^*$ to the Lie algebra of $G$;
\item[(b)] provided that $G$ satifies
the condition $\codim \Sing\g^*\ge 3$ of the argument translation
method (see Example \ref{20.70}), one gets the following complete
involutive with respect to any fixed $\e^{t_0}, t_0\not\in E$, set
of functions on $M$:
\[
\mathcal{G}^{t_0}:=\sum_{s\in\mathbb{K}^2\setminus
E}\mu_s^*(Z^{\e_{\mathrm{can}}})+
\mu_{t_0}^*(Z^{\{\e^t_{\mathrm{AT}}\}}).
\]
\end{enumerate}
\end{coro}

\begin{pf*}{Proof} of (a) follows from Corollary \ref{10.120} and
from the definition of the first integrals;  proof of (b) is a
consequence of (a) and Proposition \ref{10.150}. \qed\end{pf*}
\abz\label{30.exa}
\begin{exa}\rm
Let $M=\R^{2N}$ with coordinates $\{p_j,q_j\}_{j=1}^N$,
$\e_1=\frac{\d}{\d p_1}\wedge\frac{\d}{\d q_1}+\cdots+\frac{\d}{\d
p_N}\wedge\frac{\d}{\d q_N}$, $\e_2=a_1\frac{\d}{\d
p_1}\wedge\frac{\d}{\d q_1}+\cdots+a_N\frac{\d}{\d
p_N}\wedge\frac{\d}{\d q_N}$, where $a_1,\ldots,a_N$ are different
real numbers. Then the family $\{\e^t\}:=\{t_1\e_1+t_2\e_2\},
t=(t_1,t_2)\in \R^2$, is a dull micro-Jordan bi-Poisson structure
with the exceptional set $E=\langle
(a_1,-1)\rangle\cup\cdots\cup\langle (a_N,-1)\rangle\subset\C^2$,
the exceptional bivectors
$\widetilde{\e}^{e_j}=\widetilde{\e}^{(a_j,-1)}$ and the
corresponding foliations of symplectic leaves
$\mathcal{S}_j=\{P_j={\mathrm{const}},Q_j={\mathrm{const}}\},
j=1,\ldots,N$, where $\{P_j,Q_j\}_{j=1}^N, P_j=p_j+\i \hat{p}_j,
Q_j=q_j+\i \hat{q}_j$, are the holomorphic coordinates on
$M^\C=\C^{2N}$.

Assume $G=\mathrm{SL(2,\R)}$ is acting on $\R^2$ in a standard
linear way and that this action is extended to $M=\R^{2N}$
diagonally. It is easy to see that all these data satisfy the
assumptions of Theorem \ref{main1}. Moreover, the stabilizers
$\widetilde{G}_1,\ldots,\widetilde{G}_N\subset\widetilde{G}=\mathrm{SL(2,\C)}$
of fixed symplectic leaves $S_j=\{P_j=b_j,Q_j=c_j\}\subset\S_j,
j=1,\ldots,N$, which coincide
with the stabilizers of the vectors $\left[\begin{array}{c} b_j \\
c_j\end{array}\right]$ under the standard linear
$\widetilde{G}$-action, are 1-dimensional, consequently abelian
and have rank 1 equal to rank of $\widetilde{G}$. Hence the
reduced bi-Poisson structure $\{(\e^t)'\}$ is Kronecker on the
regular part of the variety $M/G$.

The calculations show that the moment map which corresponds to
$\e^t$ is
\[
\mu_t:(p,q)\mapsto\left[\begin{array}{ccr}z_1&=&-\sum_j\frac{p_jq_j}{t_1+a_jt_2}\\
z_2&=&-(1/2)\sum_j\frac{q_j^2}{t_1+a_jt_2}\\
z_3&=&(1/2)\sum_j\frac{p_j^2}{t_1+a_jt_2}
\end{array}\right]:\R^{2N}\rightarrow (\mathrm{sl(2,\R)})^*
\]
and that the Casimir function of $\e^{\mathrm{can}}$ on
$(\mathrm{sl(2,\R)})^*$ is $f=z_1^2+4z_2z_3$. Introducing the
affine parameter $r=-(t_1/t_2)$ we get an involutive family of
functions on $M$:
\[
p^*\mathcal{F}=\sum_{r\in
\R}\left<(\sum_{j=1}^N\frac{p_jq_j}{r-a_j})^2-(\sum_{j=1}^N\frac{q_j^2}{r-a_j})(\sum_{j=1}^N\frac{p_j^2}{r-a_j})\right>
\]
(here $p:M\rightarrow M/G$ is the canonical map). Expanding this
expression with respect to the powers of $r-a_j$ and calculating
the coefficients corresponding  to the first powers we obtain the
following functions generating $p^*\mathcal{F}$:
\[
\sum_{k=1,k\not =j}^N\frac{(p_kq_j-p_jq_k)^2}{a_k-a_j},
j=1,\ldots, N.
\]
There is one relation between these functions. By Corollary
\ref{30.cor} (b) applied with the choice $t_0=(1,0)$ (i.e.
$\e^{t_0}=\e_1$)  $p^*\mathcal{F}$ can be completed by the
function $\mu_{(1,0)}^*g$, where $g=z_1z_1^0+2z_2z_3^0+2z_3z_2^0$
is obtained from $f$ by the shift in the direction of an element
$z^0=(z_1^0,z_2^0,z_3^0)\in (\mathrm{sl(2,\R)})^*$:
\[
f(z+\la z^0)=f(z)+2\la g(z)+\la^2f(z^0), \la\in\R.
\]
Finally, we get the following complete involutive (with respect to
a standard Poisson bracket) set of functions on $\R^{2N}$:
\[
\sum_{k=1,k\not =j}^N\frac{(p_kq_j-p_jq_k)^2}{a_k-a_j},
j=1,\ldots, N-1, z_1^0\sum_{j=1}^N p_jq_j+z_2^0\sum_{j=1}^N
p_j^2+z_3^0\sum_{j=1}^N q_j^2,
\]
where $z_i^0,i=1,2,3$, are any constants simultaneously not equal
to 0.
\end{exa}

\section{Main example: diagonal action of a Lie
group on the product of $N$ copies of the dual space to its Lie
algebra} \label{s40}
\abz\label{40.10} Let $G$ be a complex Lie group, $\g$ its Lie
algebra. There is a natural coadjoint action of the direct product
$G^{\times N}$ of $N$ copies of $G$ on $(\g^*)^{\times N}$ which
restricts to $G\subset G^{\times N}$ embedded diagonally. Let
$p_j:(\g^*)^{\times N}\rightarrow \D, j=1,\ldots ,N$, denote the
natural projection to the j-th component and let $\e$ be the
canonical linear Poisson bivector (c.l.P.b) on $\D$. Then the
c.l.P.b. $\e^{\times N}$ on $(\g^*)^{\times N}$ has the
decomposition $\e^{\times N}=\e_{(1)}+\cdots+\e_{(N)}$, where
$\e_{(j)}, j=1,\ldots ,N$, is the unique Poisson bivector on
$(\g^*)^{\times N}$ defined by the condition $p_{j*}\e_{(j)}=\e,
p_{i*}\e_{(j)}=0, i\not=j$.
\abz\label{40.20}
\begin{propo}
Fix a coadjoint orbit $\O=G^{\times N}(x_1,\ldots
,x_N)\subset (\g^*)^{\times N}$ of an element $(x_1,\ldots
,x_N)\in(\g^*)^{\times N}$ and different numbers
$a_1,\ldots ,a_N\in \C$.  Then
\begin{enumerate}
\item[(a)]
the bivectors $\e^{\times N}$ and $\e^{a\times
N}=a_1\e_{(1)}+\cdots+a_N\e_{(N)}$ form a Poisson pair on
$(\g^*)^{\times N}$;
\item[(b)] they are $G$-invariant;
\item[(c)] they have the natural
restrictions (being Poisson bivectors) $\e_1=\e^{\times
N}|_\O,\lbr \e_2=\e^{a\times N}|_\O$ to $\O$; \item[(d)] the
family $\{\e^t=t_1\e_1+t_2\e_2\}$ is a dull micro-Jordan
bi-Poisson structure on $\O$ with the exceptional set $E=\langle
(a_1,-1)\rangle\cup\cdots\cup\langle (a_N,-1)\rangle$; \item[(e)]
for any $j=1,\ldots,N$ the symplectic foliation $\S_j$ of the
exceptional bivector $\e^{e_j}=\e^{(a_j,-1)}$ coincides with the
foliation of fibers of the natural projection
$p_j|_\O:\O=Gx_1\times\cdots\times Gx_N\rightarrow Gx_j$.
\end{enumerate}
\end{propo}

\begin{pf} Item (a) follows from Proposition \ref{10.45} since
$[\e_{(i)},\e_{(j)}]=0$ for any $i,j=1,\ldots ,N$.
The first bivector is $G$-invariant by definition. The
invariance of the second one follows from the
$G$-equivariance of the projections $p_j$ and from the
invariance of $\e$.

The restriction of
$\e^{\times N}$ to $\O$ is simply the restriction to a
symplectic leaf.  Moreover, any $\e_{(j)}$ is tangent to
the leaves of any projection $p_i, i\not=j$, and to
$p_j^{-1}(Gx_j)$, i.e.  $\e_{(j)}$ also has the restriction
to $\O$. This implies (c).

Since $\e_1$ is nondegenerate (as any restriction of a Poisson
bivector to a symplectic leaf), $\{\e^t\}$ is micro-Jordan.
Obviously, the only degenerate bivectors in this family are those
proportional to $\e^{e_j},j=1,\ldots ,N$, and the cooresponding
characteristic distributions satisfy the equalities
$\chi^{\e^{e_j}}=\sum_{i\not=j}(\chi^{\e_{(i)}})|_\O$, which
complete the proof. \qed\end{pf}

The main result of this section (Theorem \ref{main2}) will study
the reduction of the bi-Poisson structure $\{\e^t\}$ on $G^{\times
N}$-orbits uder the action of $G$. Now we shall
specify the class of orbits under consideration.

\abz\label{adm}
\begin{defi}
An orbit $\O=G^{\times N}(x_1\ldots ,x_N)\subset (\g^*)^{\times
N}$ is called admissible if: \begin{enumerate} \item there exist
elements $x_1'\in Gx_1,\ldots,x_N'\in Gx_N$ such that their
stabilizers $G^{x_j'}_{\g^*}\subset G, j=1,\ldots,N$, have
discrete intersection; equivalently:
\[
\g^{x_1'}_{\g^*}\cap\cdots\cap \g^{x_N'}_{\g^*}=\{0\};
\]
\item the stabilizers $G_j:=G^{x_j}_{\g^*}\subset G, j=1,\ldots,N$,
have all the same rank equal to the rank of $G$:
\[
\rk G_1=\cdots=\rk G_N=\rk G.
\]
\end{enumerate}
\end{defi}

We postpone the discussion of the question which orbits are
admissible to the end of this section (see Subsection
\ref{dissc}); here we mention only that the admissibility holds
for generic orbits in the semisimple case.

Now we formulate the second main result of this paper.
\abz\label{main2}
\begin{theo}
Let $\O\subset (\g^*)^{\times N}$ be an admissible $G^{\times
N}$-orbit and let $M\subset \O$ be an open set such that $M/G$ is
a manifold.  Then the bi-Poisson structure $\{\e^t\}|_M$ is
projectable via the canonical map $p:M\rightarrow M/G$ and the
projection $\{(\e^t)'\}$ is a micro-Kronecker bi-Poisson structure
(see Definition \ref{20.50}) on $M'=M/G$. More precisely,
$\{(\e^t)'\}$ is Kronecker at any $x'\in M'\setminus
p(\mathcal{N})$, where $\mathcal{N}\subset(\g^*)^{\times N}$ is
the algebraic set of all elements with a nondiscrete
$G$-stabilizer.
\end{theo}

\begin{pf} Of course, this proof will use
Theorem \ref{main1}. Now we shall check that the
$G$-action on $\{\e^t\}$ satisfies the assumptions of this
theorem.

First, we note that since the $G$-stabilizer
$G^{(x_1,\ldots,x_N)}_{(\g^*)^{\times N}}$  of a point
$(x_1,\ldots,x_N)\in (\g^*)^{\times N}$ is equal to the
intersection $G^{x_1}_{\g^*}\cap\cdots\cap G^{x_N}_{\g^*}$,
condition (1) in definition of admissibility guarantees that the
$G$-action is locally free.

 The $G$-invariance of $\{\e^t\}$ was proved
in \ref{40.20}(a), so we get assumption (1) of Theorem
\ref{main1}. To check assumption (2) recall (see \ref{40.20}(e))
that the symplectic foliation of the exceptional bivector
$\e^{e_j}$ coincides with $\{Gx_1\times\cdots \times
Gx_{j-1}\times x\times Gx_{j+1}\times\cdots\times Gx_N\mid x\in
Gx_j\}$. Since $G$ is acting transitively on $Gx_j$, the same is
true for the induced $G$-action on the leaves of this foliation.

Now, let us prove the hamiltonicity of the $G$-action on
$M$ with respect to $\e^t, t\in\C^2\setminus E$. The
commutativity of the following diagram is standard:
\[
\begin{array}{ccc}
\g^{\times N} & \stackrel{i}{\longrightarrow}&
\E((\g^*)^{\times N}) \\ \parallel & & \downarrow
\e^{\times N}(\cdot) \\ \g^{\times N}
&\stackrel{\rho}{\longrightarrow} & \G T((\g^*)^{\times
N})
\end{array}
\]
(here $i$ is  the inclusion of $\g^{\times N}$ in
$\E((\g^*)^{\times N})$ as a set of linear functions and
$\rho$ is the Lie algebra homomorphism corresponding to the
coadjoint action). It leads to the following commutative
diagram:
\[
\begin{array}{ccc}
\g & \stackrel{\psi^t}{\longrightarrow}& \E(\O) \\ \parallel & &
\downarrow \e^t(\cdot) \\ \g &\stackrel{\rho^d}{\longrightarrow} &
\G T\O,
\end{array}
\]
where $\rho^d$ is the restriction of $\rho$ to the diagonal,
$\psi^t$ is defined as
\[
\psi^t(x)=\frac{1}{t_1+a_1t_2}p_1^*(x)|_\O+
\cdots+\frac{1}{t_1+a_Nt_2}p_N^*(x)|_\O,
\]
$x$ in the RHS being understood as a function on $\D$. So
assumption \ref{main1}(4) is satisfied, it remains to check
\ref{main1}(3). This will be done with the help of the commutative
diagram
\[
\begin{array}{ccc}
\g_j & \stackrel{\psi_j}{\longrightarrow}&
\E(S_j) \\ \parallel & & \downarrow
\e^{e_j}(\cdot) \\ \g_j
&\stackrel{\rho^d|_{\g_j}}{\longrightarrow} & \G T\S_j.
\end{array}
\]
Here $\g_j$ is the Lie algebra of the stabilizer
$G_j=G_{\g^*}^{x_j}$ of a symplectic leaf $S_j=Gx_1\times\cdots
\times Gx_{j-1}\times x_j\times Gx_{j+1}\times\cdots\times Gx_N$,
$\rho^d|_{\g_j}$ is the restriction to $\g_j$ of the above
mentioned map $\rho^d$, and $\psi_j$ is given by the formula
\begin{eqnarray*}
\psi_j(x)&=&\frac{1}{a_j-a_1}p_1^*(x)|_{S_j}+\cdots+
\frac{1}{a_j-a_{j-1}}p_{j-1}^*(x)|_{S_j}\\ & &
+\frac{1}{a_j-a_{j+1}}p_{j+1}^*(x)|_{S_j}+
\cdots+\frac{1}{a_j-a_N}p_N^*(x)|_{S_j},\ x\in\g_j\subset\E(\g^*).
\end{eqnarray*}

Thus, all the assumptions of Theorem \ref{main1} are checked. In
order to finish the proof we need to use condition (2) of equality
of ranks from the definition of admissibility. \qed\end{pf}

\abz\label{40.60}
\begin{coro}
The moment map $\mu_t:\O\longrightarrow\g^*$ for the $G$-action on
$(\O,\e^t)$ is given by the restriction to $\O$ of the following
map:
\[
(\g^*)^{\times N}\ni (x_1,\ldots
,x_N)\mapsto\frac{1}{t_1+a_1t_2}x_1+\cdots+\frac{1}{t_1+a_Nt_2}x_N.
\]
\end{coro}
\begin{pf*}{Proof} follows from the proof of Theorem \ref{main2}. \qed\end{pf*}

\abz\label{40.70}
\begin{coro}
The set of first integrals $Z^{\{(\e^t)'\}}$ of the reduced
Kronecker bi-Poisson structure coincides with the family of
functions
\[
\mathcal{F}=\sum_{t\in\C^2\setminus
E}\mu_t^*(Z^{\e_{\mathrm{can}}})
\]
considered as functions on $M/G$.
\end{coro}

\begin{pf*}{Proof} follows from Corollary \ref{10.120}. See also Corollary
\ref{30.cor}. \qed\end{pf*}
\abz\label{40.80}
\begin{coro}
Assume that $G$ satifies the condition $\codim \Sing\g^*\ge 3$ of
the argument translation method (see Example \ref{20.70}). Then
for any fixed $t_0\in \C^2\setminus E$ and any regular $a\in \g^*$
we get a complete involutive set of functions on $\O$
\[
\mathcal{G}^{t_0}=\sum_{t\in\C^2\setminus
E}\mu_t^*(Z^{\e_{\mathrm{can}}})+\mu_{t_0}^*(Z^{\{\e^t_{\mathrm{AT}}\}}).
\]
\end{coro}

\begin{pf*}{Proof} follows from Corollary \ref{30.cor}. \qed\end{pf*}

In the remaining part of this section we want to discuss two
aspects of applicability of Theorem \ref{main2}, i.e. which orbits
are admissible and what happens in real case.
\abz\label{dissc}
\begin{theo}
Assume $G$ is semisimple. Then a generic $G^{\times N}$-orbit
$\mathcal{O}=Gx_1\times\cdots\times Gx_N\subset (\g^*)^{\times N}$
is admissible for any $N\ge 2$.
\end{theo}

\begin{pf} We will first prove condition (2) of Definition \ref{adm}. It
follows from the well known fact (see \cite{ag} for example), that
the stabilizers of generic elements in the dual space to any Lie
algebra are abelian, and from the equality of dimensions: $\rk
G=\dim G_{\D}^{x_1}=\cdots=\dim G_{\D}^{x_N}$.

The first condition of the definition of admissibility requires
some additional preparations.
\end{pf}
\abz\label{40.40}
\begin{lemm}
Let $K\subset G$ be a maximal compact subgroup. Then the principal
orbital type stabilizer $K_{\D}^x\subset K$ of an element $x\in
\g^*$ under the coadjoint action of $K$ on $\D$ is at most
discrete (finite).  \end{lemm}

\begin{pf} (The idea of this proof was communicated to the author by
Prof. Sam Evans.) For this proof we identify $\D$ and $\g$ using
the Killing form. We claim that the Lie algebra $\k^x$ of a
principal orbital type stabilizer $K^x$ for the $K$-action on $\D$
is trivial.  Indeed, Theorem 3.6 of \cite{kost} shows that for any
nilpotent element $e\in\g$ the subalgebra $\g_e=\ad e(\g)\bigcap
\g^e$ consists of nilpotent elements.  If, moreover, $e$ is a
principal nilpotent element (see \cite{kost}, Subsection 5.2) it
can be easily seen that $\g_e=\g^e$. However, each element of $\k$
is semisimple; thus $\k^e=\k\cap\g^e=\{0\}$. Of course, this
implies the triviality of $\k^x$. \qed\end{pf}

\begin{pf*}{Continuation of the proof}
Now we are able to complete the proof of Theorem \ref{dissc}.
Since the $K$-action on $\D=\k\oplus \i\k$ is diagonal, it follows
from the above lemma that for a generic pair
$(a,b)\in\k^*\oplus\k^*$ the intersection of stabilizers
$K_{\k^*}^a\cap K_{\k^*}^b$ is finite. The complexification gives
the discreteness of the intersection $G_{\D}^a\cap G_{\D}^b$ for a
generic pair $(a,b)\in \D\oplus\D$. This implies the result.
\qed\end{pf*}

\abz\label{40.90}
\begin{rema}\rm
Theorem \ref{dissc} shows that Theorem \ref{main2} can be applied
to semisimple Lie groups and generic orbits in $(\g^*)^{\times
N}$. We also note that:
\begin{enumerate}
\item Corollary \ref{40.80}
is also valid for this data since the condition
$\codim\Sing\g^*\ge 3$ of the argument translation method (see
Example \ref{20.70}) holds in semisimple case;
\item Theorem \ref{main2} can be also applied  for nonsemisimple
Lie groups: condition (2) of definition of admissibility \ref{adm}
holds for any Lie algebra $\g$ and for the stabilizers $G^{x_j}$
of generic points $x_j\in \g^*$ (see proof of Theorem
\ref{dissc}); condition (1) should be achieved at least for the
algebras with the trivial center by increasing the number of
components $N$;
\item Another possibility for application of Theorem \ref{main2}
are nongeneric orbits: for example, rank of the stabilizer $G^x$
of any semisimple element $x\in\D$ coincides with $\rk G$ for
semisimple $G$ (see \cite{cmg}, Chapter 2, for example).
\end{enumerate}
\end{rema}

\abz\label{40.95}
\begin{rema}\rm Since the complexification of a real
semisimple Lie group is complex semisimple, all the results of
this section are valid in real setting, i.e. for a real semisimple
group $G$ and different $a_1,\ldots ,a_N\in \R$. All proofs remain
the same, only the arguments concerning the proof of condition (1)
of the definition of the  admissibility for generic orbits require
additional considerations.
\end{rema}
\abz\label{40.100}
\begin{propo}
 Let $G$ be a real semisimple Lie group with
the Lie algebra $\g$. Then the generic stabilizer of the
$G$-action on $(\g^*)^{\times N},N\ge 2$, is at most discrete.
\end{propo}

\begin{pf} Let $\g^\C$ be the complexification of $\g$. Then by Lemma
\ref{40.40} the set $\mathcal{N}$ of all points
$x\in(\g^\C)^*\times(\g^\C)^*$ with the nontrivial stabilizer
$(\g^\C)^x$ (with respect to the diagonal action of $\g^\C$) is a
proper complex algebraic set. The intersection
$\mathcal{N}'=\mathcal{N}\cap\g^*\times\g^*$ is a proper real
algebraic set, and for $x\in\g^*\times\g^*\setminus\mathcal{N}'$
the corresponding real stabilizer $\g^x=(\g^\C)^x\cap\g$ is
trivial. \qed\end{pf}
\begin{ack}
The author wishes to thank Prof. Ilya Zakharevich for critical
remarks on the previous version of this paper \cite{p1}, which
inspired the appearance of the present one. The author is also
indebted to Prof. Sam Evans for the idea of the proof of
Proposition \ref{40.40} and to Prof. Pantelis Damianou for
indicating the reference \cite{cmg}.
\end{ack}


\begin{thebibliography}{11}

\bibitem{ag}
V. Arnold and A. Givental, Symplectic geometry, in: {\em
Encyclopaedia of Mathematical Sciences (Dynamical systems
  IV), volume~4\/} (Springer, 1990) 1--136.

\bibitem{bols}
A. Bolsinov. Compatible poisson brackets on Lie algebras and
completeness of
  families of functions in involution,
{\em Izv.Akad. Nauk SSSR, Ser. Mat. \/} {\bf 55} (1991), in
Russian. English translation: {\em Math. USSR Izvestiya\/} {\bf
38} (1992) 69--90.

\bibitem{bw}
F. Bruhat and H. Whitney, Quelques propri\'{e}t\'{e}s
fondamentales des ensembles
  analitiques-r\'{e}els,
{\em Comment.Math.Helv.\/} {\bf 33} (1959) 132--160.

\bibitem{cmg}
D. Collingwood and W. McGovern, {\em Nilpotent orbits in
semisimple Lie algebras\/} (VNR, 1993).

\bibitem{g1}
M. Gaudin, Mod\`{e}les exacts en m\'{e}canique statistique: la
m\'{e}thode de
  Bethe et ses g\'{e}n\'{e}ralisations,
{\em Note CEA} {\bf 1559 (1)} et {\bf 1559 (2)} (1972).

\bibitem{g2}
M. Gaudin, Diagonalisation d'une class d'hamiltoniens de spin,
{\em J. de Physique\/} {\bf 37} (1976) 1087--1098.

\bibitem{g3}
M. Gaudin, {\em La fonction d'onde de Bethe\/} (Paris: Masson,
1983).

\bibitem{gz1}
I. Gelfand and I. Zakharevich, Spectral theory for a pair of
skew-symmetrical operators on $S^1$, {\em Funktsion. Analiz i ego
Prilozh.\/} {\bf 23 (2)} (1989) 1--11, in Russian. English
translation: {\em Functional Anal. Appl.\/} {\bf 23(2)} (1989)
85--93.

\bibitem{gz2}
I. Gelfand and I. Zakharevich, Webs, Veronese curves, and
bihamiltonian systems, {\em J. Funkt. Anal.\/} {\bf 99} (1991)
150--178.

\bibitem{gz3}
I. Gelfand and I. Zakharevich, Webs, Lenard schemes, and the local
geometry of bihamiltonian Toda
  and Lax structures.
{\em Selecta-Math. (N.S.)\/} {\bf 6} (2000) 131--183.

\bibitem{gs}
V. Guillemin and S. Sternberg, Convexity properties of the moment
map, {\em Invent. Math.\/} {\bf 67} (1982) 491--513.

\bibitem{kir}
A. Kirillov, Local Lie algebras, {\em Uspekhi Mat. Nauk\/} {\bf
31} (1976) 57--76, in Russian. English translation: {\em Russian
Math. Surveys\/} {\bf 31} (1976) 55--75.

\bibitem{kost}
B. Kostant, The principal three-dimensional subgroup and the Betti
numbers of a
  complex simple Lie group.
{\em Amer.J.Math.\/} {\bf 81} (1959) 973--1032.

\bibitem{lieb}
P. Libermann, Sous-varietes et feulletages symplectiquement
reguliers, in: {\em Symplectic geometry (Tolouse 1981)}, volume~80
of {\em Res.
  Notes in Math.} (Pitman, 1983) 81--106.

\bibitem{lichn}
A. Lichnerowicz, Les vari\'{e}t\'{e}s de Jacobi et leurs
alg\`{e}bres de Lie
  associ\'{e}es,
{\em J. Math. Pures Appl.\/} {\bf 57} (1978) 453--488.

\bibitem{p}
A. Panasyuk, Veronese webs for bihamiltonian structures of higher
corank, {\em Banach Center Publications\/} {\bf 51} (2000)
251--261.

\bibitem{p1}
A.Panasyuk, Symplectic realizations of bihamiltonian structures,
{\em Archieved as DG/0001126}.

\bibitem{rst}
A. Reyman and M. Semenov-Tian-Shansky, Group-theoretical methods
in the theory of integrable systems, in: {\em Encyclopaedia of
Math. Sciences (Dynamical Systems VII),
  volume~16\/} (Springer, 1994)  116--225.

\bibitem{s1}
E. Sklyanin, Separation of variables in the Gaudin model, {\em
Zap. Nauchn. Sem. Leningrad. Otdel. Math. Inst. Steklov.
  (LOMI)\/} {\bf 164} (1987) 151--169, in Russian. English translation:
  {\em J. Soviet Math.\/} {\bf 47} (1987) 2473--2488.

\bibitem{s2}
E. Sklyanin, Separation of variables in the classical integrable
$\mathrm{sl(3)}$ magnetic
  chain,
{\em preprint RIMS-871} (1992).

\bibitem{s3}
E. Sklyanin, Separation of variables -- new trends, {\em Progr.
Theoret. Phys. Suppl.\/} {\bf 118} (1995) 35--60.

\bibitem{weinst}
A. Weinstein, The local structure of Poisson manifolds, {\em J.
Diff. Geom.\/} {\bf 18} (1983) 523--557.

\bibitem{z1}
I. Zakharevich, Kronecker webs, bihamiltonian structures, and the
method of argument
  translation,
{\em Transformation Groups\/} {\bf 6} (2001) 267--300.

\end{thebibliography}
\end{document}